\documentclass[12pt]{article}

\usepackage[english]{babel} 
\usepackage{graphicx}
\usepackage{amsmath, amssymb, amsthm, mathtools, stmaryrd}

\usepackage[letterpaper, margin=1in]{geometry}

\usepackage{graphicx, float, subcaption, tikz}
\usetikzlibrary{decorations.pathreplacing, decorations.markings, arrows.meta, positioning}

\usepackage[numbers]{natbib} % Use numeric citations
\usepackage{enumitem, stackengine, fancyhdr}

\usepackage[colorlinks=true, allcolors=blue]{hyperref}

\theoremstyle{definition}
\newtheorem{definition}{Definition}[section]

\theoremstyle{definition}
\newtheorem{theorem}[definition]{Theorem}
\theoremstyle{definition} 
\newtheorem{lemma}[definition]{Lemma}
\theoremstyle{definition} 
\newtheorem{corollary}[definition]{Corollary}
\theoremstyle{definition} 

\theoremstyle{definition} 
\newtheorem{example}[definition]{Example}
\theoremstyle{definition} 
\newtheorem{remark}[definition]{Remark}

\newcommand{\R}{\mathbb{R}}
\newcommand{\Q}{\mathbb{Q}}
\newcommand{\N}{\mathbb{N}}
\newcommand{\Z}{\mathbb{Z}}
\newcommand{\C}{\mathbb{C}}

\title{Certain Representation of \(\Q\) With Bounded Partial Quotients}\author{Franciszek Kobus \\ \small{UCL, Department of Mathematics}}
\date{\today}

\begin{document}
\maketitle

\begin{abstract}
We define Collatz representations for a subset of rational numbers and prove that each real number \(x \notin (-1,1)\) can be approximated arbitrarily well by rational numbers which have only \(2\)'s and \(1\)'s in their Collatz representation.
\end{abstract}

\section{Introduction}

In this short article we define an original way of representing rational numbers with odd numerators and denominators \footnote{It could easily be extended to \(\mathbb{Z}\) localised at \((2)\) but for the sake of clarity we stick to the set of rationals with odd numerators and denominators.} analogous to continued fractions of rational numbers, and continued fractions of polynomials (see \cite{Malagoli}). The definition is of course also motivated by the famous Collatz conjecture (\cite{Collatz}), however this article relies more on its connection to the other definitions of continued fractions in different settings. As pointed out by Marshall Hall in \cite{Hall1947OnTS}, restricting our attention to the set of continued fractions of rationals with bounded partial quotients yields quite unexpected results. In our setting we focus on the set of rationals whose Collatz representations consist only of \(1\)'s and \(2\)'s. We prove that this set is dense in \((-\infty,-1]\cup [1,\infty)\) and exhibits a particular fractal structure. At the end of the article we list some motivating further results.

\section{Collatz representation}
\definition By \(\Q^{odd}\) we denote the set of rational numbers with odd numerators and odd denominators. By \(S^{inf}\) we denote the set of all infinite sequences of positive integers. \(S^{fin}\) is the set of all finite sequences of positive integers.
\definition The \textbf{Collatz algorithm} runs as follows:

\begin{enumerate}
    \item Start with a number \(a \in \Q^{odd}\). Take \(i\coloneqq1\), \(A_1 \coloneqq 3a+1\).
    \item Denote by \(a_i\) the \(2\)-adic valuation of the numerator of \(A_i\) (written as a fraction in its reduced form). Let \(B_i \coloneqq A_i \cdot 2^{-a_i}\).
    \item Take \(i\coloneqq i+1\), \(A_i \coloneqq 3B_{i-1}+1\). Return to step \(2\).
\end{enumerate}

\definition The \textbf{Collatz representation} of a number \(a\in \Q^{odd}\) is an infinite sequence \(C(a)=[a_1,a_2,...]\) of integers \(a_i\) obtained by running the Collatz algorithm on \(a\). We call the numbers \(a_i\) the \textbf{partial quotients} of \(a\).
\remark All \(a_i\)'s are positive integers. Observe that if \(a\in \Q^{odd}\) is of the form \(x/y\), for \(x\), \(y\) odd integers, then the numerator of \(3a+1\) in its reduced form equals \(3x+y\) if \(3\) does not divide \(y\), or \(x+y/3\) if \(3\) divides \(y\). In both of these cases it is even. On the other hand, the denominator of \(3a+1\) equals \(y\) or \(y/3\), again depending on whether or not \(y\) is divisible by \(3\). In both cases however this denominator is odd. Hence \(a_1\) is a positive integer and \(B_1 \in \Q^{odd}\). By induction, all \(a_i\)'s are positive integers.\\

Hence we have a well-defined function:
\begin{align}
    C: \Q^{odd} \rightarrow S^{inf}
\end{align}
This function will be the focus of our research. We will find the inverse of it for sequences in \(S^{inf}\) which are eventually periodic.\\

\example Running the Collatz algorithm for \(1\) yields:
\( 
1 \xrightarrow{ \frac{ 3x+1}{2^2}} 1\xrightarrow{ \frac{ 3x+1}{2^2}} 1  \dots\)
Hence \(C(1)=[(2)]\).
\example Running the Collatz \(
-1 \xrightarrow{ \frac{ 3x+1}{2^1}} -1\xrightarrow{ \frac{ 3x+1}{2^1}} -1 \dots\)
Hence \(C(-1)=[(1)]\).\\

Thanks to these examples we can observe that if a number \(a\in \Q^{odd}\) eventually reaches \(1\) in its Collatz algorithm after a finite number of steps, this means that there is \(N\in \N\), such that \(a_n=2, \forall n\geq N\). On the other hand, if it eventually reaches \(-1\), then there is \(N'\in \N\), such that \(a_n=1, \forall n\geq N'\). In general, we see that if there exist \(k\), \(n \in \N\) such that \(A_k=A_{n+k}\), then \(C(a)=[a_1,...,a_n,(a_{n+1},...,a_{n+k})]\). We call such numbers \(a\) \textbf{absolutely periodic} if 
\(k=0\) and \textbf{eventually periodic} if \(k>0\). The number \(n\) is called the \textbf{period} of \(a\).

\example

\(5/3 \xrightarrow{ \frac{ 3x+1}{2^1}} 3\xrightarrow{ \frac{ 3x+1}{2^1}} 5 \xrightarrow{ \frac{ 3x+1}{2^4}} 1\xrightarrow{ \frac{ 3x+1}{2^2}} 1\xrightarrow{ \frac{ 3x+1}{2^2}} 1 \dots\)
Hence \(C(5/3)=[1,1,4,(2)]\).

\example\label{examplesevenfive}
\(
-7/5 \xrightarrow{ \frac{ 3x+1}{2^4}} -1/5 \xrightarrow{ \frac{ 3x+1}{2^1}} 1/5 \xrightarrow{ \frac{ 3x+1}{2^3}} 1/5\xrightarrow{ \frac{ 3x+1}{2^3}} 1/5 \dots\)
Hence \(C(-7/5)=[4,1,(3)]\).

\example
\(-7 \xrightarrow{ \frac{ 3x+1}{2^2}} -5 \xrightarrow{ \frac{ 3x+1}{2^1}} -7 \xrightarrow{ \frac{ 3x+1}{2^2}} -5 \dots\)
Hence \(C(-7)=[(2,1)]\), and \(C(-5)=[(1,2)]\).\\

\lemma
For any eventually periodic sequence \([a_1,...,a_m,(b_1,...,b_n)] \in S^{inf}\) there exists precisely one \(a \in \Q^{odd}\) such that \(C(a)=[a_1,...,a_m,(b_1,...,b_n)]\). We denote it as \(C^{-1}([a_1,...,a_m,(b_1,...,b_n)])\).
\begin{proof}
Let's first prove the statement for absolutely periodic sequences. Take \([(b_1,...,b_n)] \in S^{inf}\). If there exists \(a \in \Q^{odd}\) such that \(C(a)=[(b_1,...,b_n)]\), then we have:
\begin{align}
    a = \underbrace{(3(...(3(3}_{\text{$n$ times}}a+1)/2^{b_1}+1)/2^{b_2}...)+1)/2^{b_n}.
\end{align}
Unfolding that, we get:
\begin{align}
    2^{b_1+...+b_n} \cdot a = 3^n \cdot a + 3^{n-1} +3^{n-2}\cdot 2^{b_1}+3^{n-3}\cdot 2^{b_1+b_2}+\dots + 2^{b_1+\dots +b_{n-1}}
\end{align}
hence:
\begin{align}
    a&= \frac{3^{n-1} +3^{n-2}\cdot 2^{b_1}+3^{n-3}\cdot 2^{b_1+b_2}+\dots + 2^{b_1+\dots +b_{n-1}}}{2^{b_1+...+b_n}- 3^n}\\
    &=\frac{3^{n-1} +(\sum_{i=1}^{n-1} (3^{n-1-i} \cdot \prod_{j=1}^i 2^{b_j}) )}{2^{b_1+...+b_n}- 3^n}\label{eqperiod}
\end{align}
So, if there is \(a\) such that \(C(a)=[(b_1,...,b_n)]\), it satisfies equation \eqref{eqperiod}, hence it's unique. It remains to show that for \(a\) as in equation \eqref{eqperiod}, when we run the Collatz algorithm on it, it actually gives us the sequence \([(b_1,...,b_n)]\). This can be easily checked by induction on the period of \(a\). A proof for the case of eventually periodic sequences runs analogously.\\
\end{proof}

\corollary The function \(C\) together with \(C^{-1}\) define a bijection between the set of all eventually periodic sequences and the set of rational numbers whose Collatz representation is eventually periodic.\\

\section{Bounded partial quotients}

The main point of this article is to show a fractal structure exhibited by the set of numbers \(C^{-1}([A])\), where \([A]\) is eventually periodic and consists only of \(1\)'s and \(2\)'s.

\definition For a finite sequence \(A=[a_1,a_2,...,a_n]\), with \(a_i \in \{1,2 \}\), \(\forall i \in \{1,2,...,n\}\) and \(n \in \Z_{\geq 0}\), we call the point \(P(A) \coloneqq (C^{-1}[A,(2)],C^{-1}[A,(1)])\in \R^2\) the point \textbf{related to} A.

\example Take A to be the empty sequence, of length zero. Then:
\[ P(A)=(C^{-1}[(2)],C^{-1}[(1)])=(1,-1).\]

\example Take A to be the sequence \([1]\), of length one. Then:
\[P(A)=(C^{-1}[1,(2)],C^{-1}[1,(1)])=(1/3,-1).\]

\example Take A to be the sequence \([2]\), of length one. Then: 
\[P(A)=(C^{-1}[2,(2)],C^{-1}[2,(1)])=(1,-5/3).\]

Let \(n\) be an integer and take all the sequences of length less or equal to $n$, with elements in \(\{1,2 \}\), and mark the points related to them on the \(\R^2\)-plane. Moreover if for some sequences \([A]\), \([B]\) we have that \([B] = [A,1]\) or \([B] = [A,2]\), connect the points \(P(A)\) and \(P(B)\) with a line segment. Call the resulting picture \(P_n\). Running a Python code for different values of \(n\) yields the following images:
\begin{figure}[H]
    \centering
    % First subfigure
    \begin{subfigure}[b]{0.45\linewidth}
        \centering
        \includegraphics[width=\linewidth]{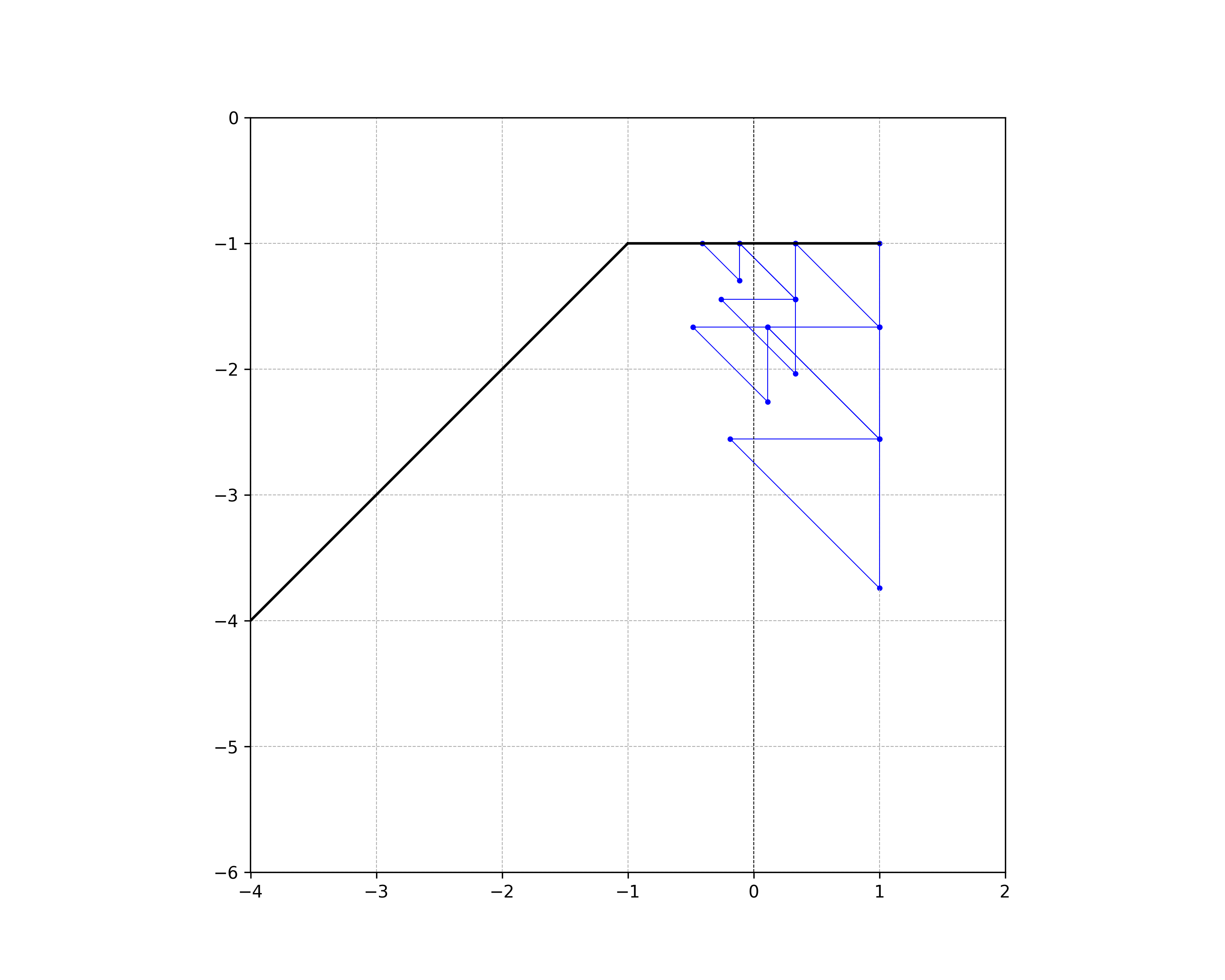}
        \caption{$P_3$}
        \label{fig:first}
    \end{subfigure}
    % Second subfigure
    \begin{subfigure}[b]{0.45\linewidth}
        \centering
        \includegraphics[width=\linewidth]{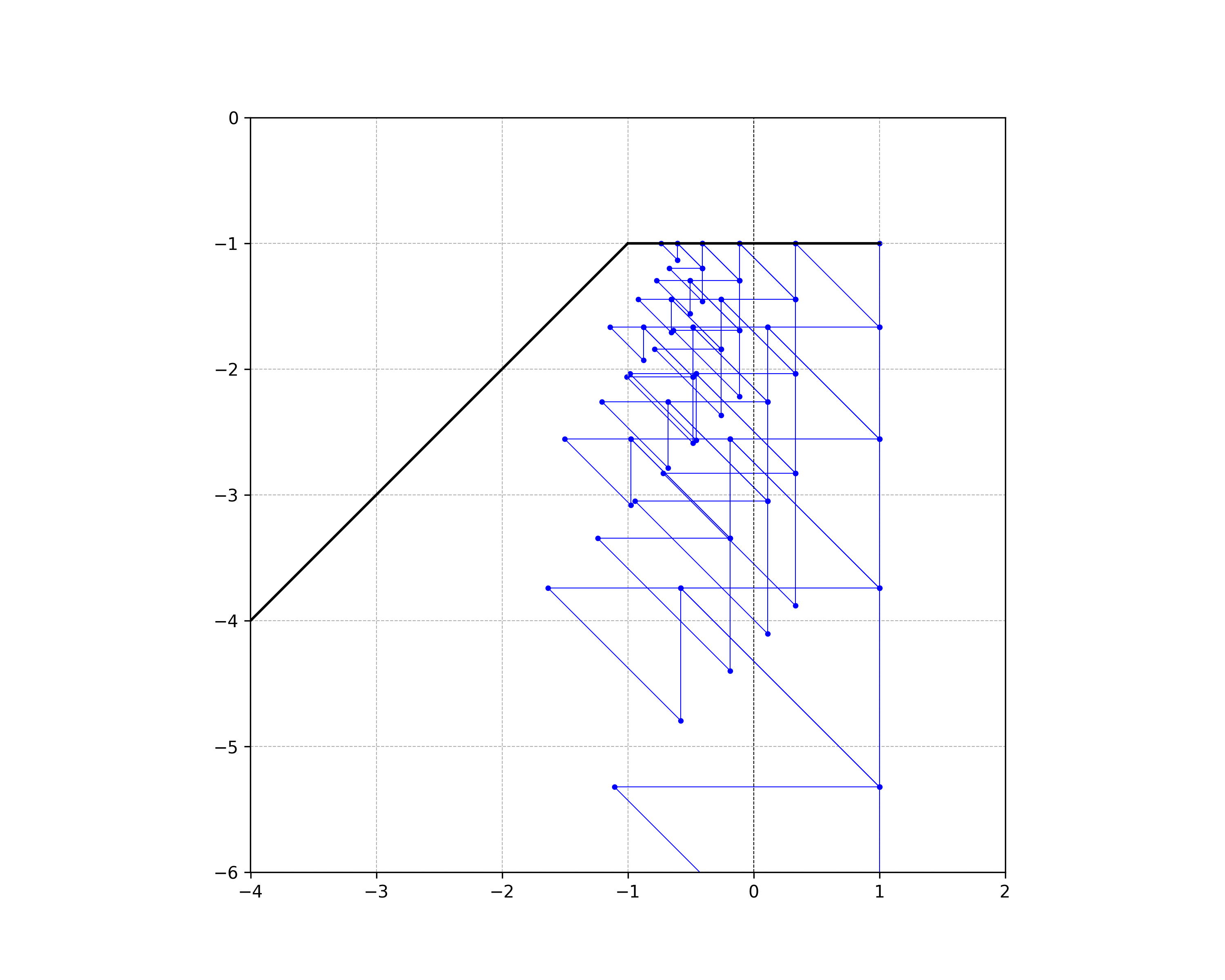}
        \caption{$P_5$}
        \label{fig:second}
    \end{subfigure}
    \label{fig:side-by-side2}
\end{figure}

\begin{figure}[H]
    \centering
     % Set the counter to 2 so the first subfigure will be (c)
    % First subfigure
    \begin{subfigure}[b]{0.45\linewidth}
    \setcounter{subfigure}{1}
        \centering
        \includegraphics[width=\linewidth]{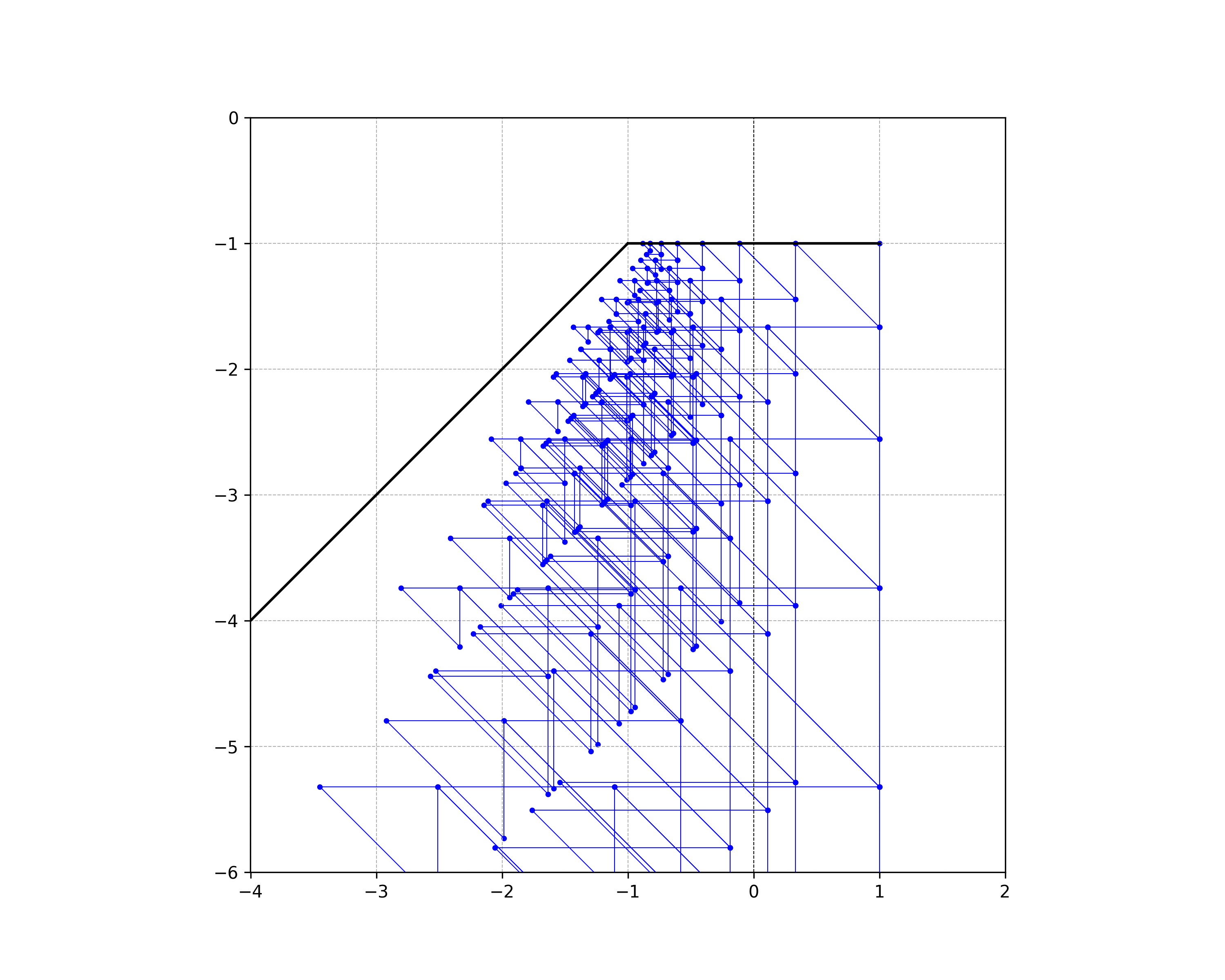}
        \caption{$P_7$}
        \label{fig:first1}
    \end{subfigure}
    \begin{subfigure}[b]{0.45\linewidth}
        \centering
        \includegraphics[width=\linewidth]{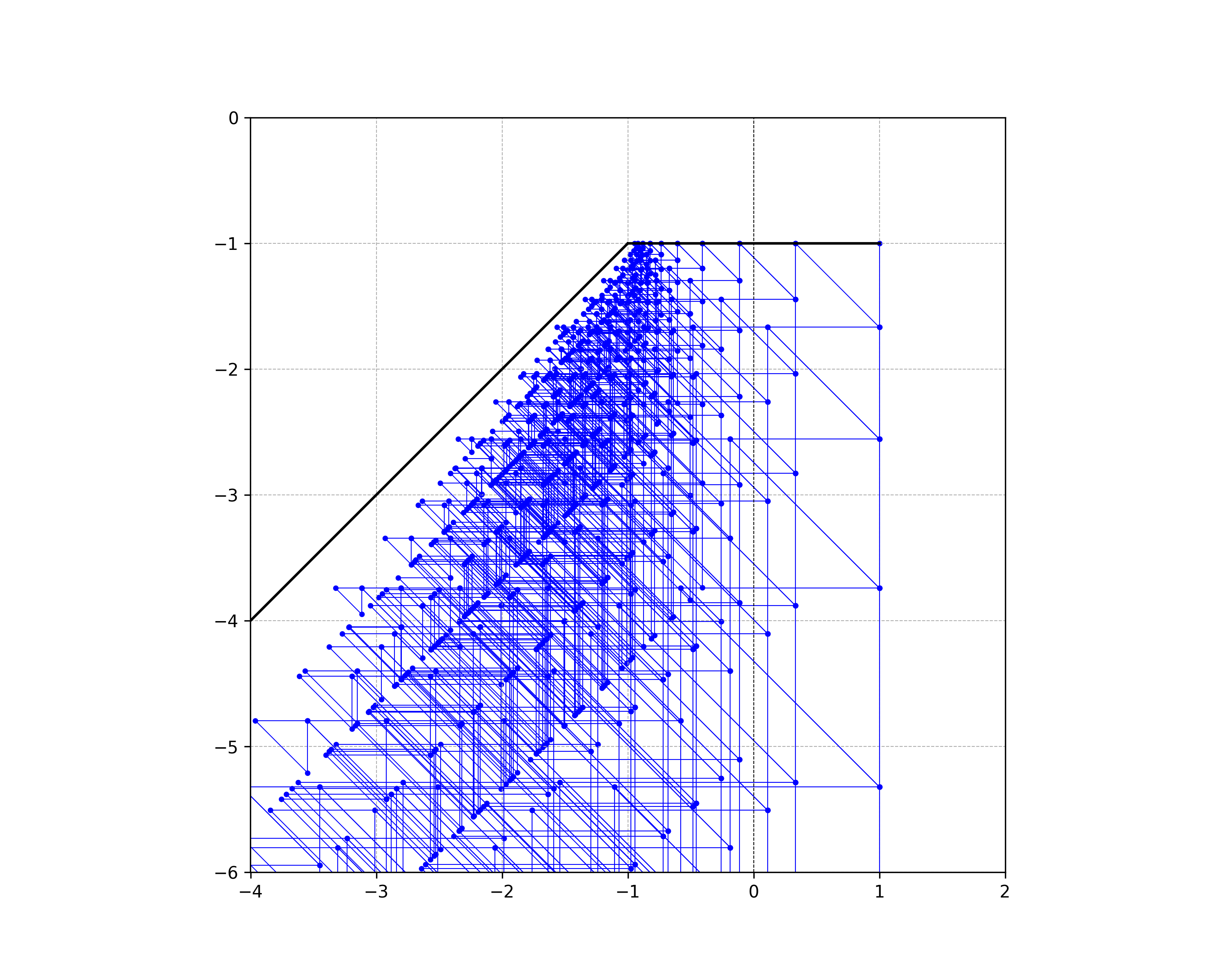}
        \caption{$P_9$}
        \label{fig:second1}
    \end{subfigure}
    
    \label{fig:side-by-side}
\end{figure}

Based on this empirical results we can informally state the following hypotheses:
\begin{enumerate}
    \item If two "triangles" we mark share a vertex, then they are similar to each other with a scaling factor \(2/3\), if they lie one next to another, and \(4/3\) if they lie one above another, and their corresponding sides are parallel to each other. 
    \item For each \(n\), \(P_n\) is enclosed above and from the right by the lines \(L_1:y-x=0\), \(L_2:x-1=0\), \(L_3:y+1=0\).
    \item The closure of \(L_1 \cap (\lim_{n \to \infty}P_n)\) in \(\R^2\) is the half-line \(L_1 \cap \{x\leq 1\}\).
\end{enumerate}

Firstly, let's prove that hypothesis \(1\). Observe that, by the uniqueness of Collatz representation, for any \(n\) we have that each point of \(P_n\) is a vertex of exactly four line segments if it's first obtained after less than \(n\) steps, and two if it's obtained precisely at the \(n\)-th step. Moreover, we have:

\lemma\label{statements} for any \([A]\) a finite sequence with partial quotients in \(\{1,2\}\) we have:
\begin{enumerate}
    \item \(C^{-1}[A,(1)]-C^{-1}[A,2,(1)]=C^{-1}[A,(2)]-C^{-1}[A,1,(2)]\).
    \item \(C^{-1}[A,2,(1)]-C^{-1}[A,2,2,(1)]= 4/3 \cdot (C^{-1}[A,(1)]-C^{-1}[A,2,(1)])\).
    \item \(C^{-1}[A,1,(2)]-C^{-1}[A,1,1,(2)]= 2/3 \cdot (C^{-1}[A,(2)]-C^{-1}[A,1,(2)])\).

\end{enumerate}
We prove these three statements using forms for the coordinates of the points of \(P_n\) analogous to \eqref{eqperiod}. Applying the statements inductively for \(n\) yields the structure described by hypothesis \(1\).\\
    
Now, hypothesis \(2\). \(P_n\) being enclosed by \(L_2\) and \(L_3\) follows directly from hypothesis \(1\). The enclosure by \(L_1\) follows from:
\theorem \(\lim_{n \to \infty}C^{-1}([A,\underbrace{1,1,...1}_{n\text{ times}}(2)]) = C^{-1}([A,(1)])\).\\

The proof relies on the fact that for a given \(n\), the difference of \(\C^{-1}([A,\underbrace{1,1,...1}_{n\text{ times}}(2)])\) and \(C^{-1}([A,(1)])\) equals
\(\frac{2}{3^m}\cdot\Big(\frac{ 2}{3}\Big)^n\), where \(m\) is the length of the sequence \(A\). This difference of course approaches \(0\) as \(n\) approaches \(\infty\).

\remark Note that it doesn't work the other way around; 
\begin{align}
    \lim_{n \to \infty}C^{-1}([A,\underbrace{2,2,...,2}_{n\text{ times}}(1)]) \neq C^{-1}([A,(2)]).
\end{align}
Even further, for any sequence \(A\), the limit on the LHS above equals \(-\infty\).\\

Let's tackle hypothesis \(3\):

\theorem Let \(z \in \R_{\leq -1}\). We then have that for every \(\epsilon>0\) there exists \(A \in S^{fin}\) consisting of only \(1\)'s and \(2\)'s, such that \(|z- C^{-1}([A,(1)]) |< \epsilon\).

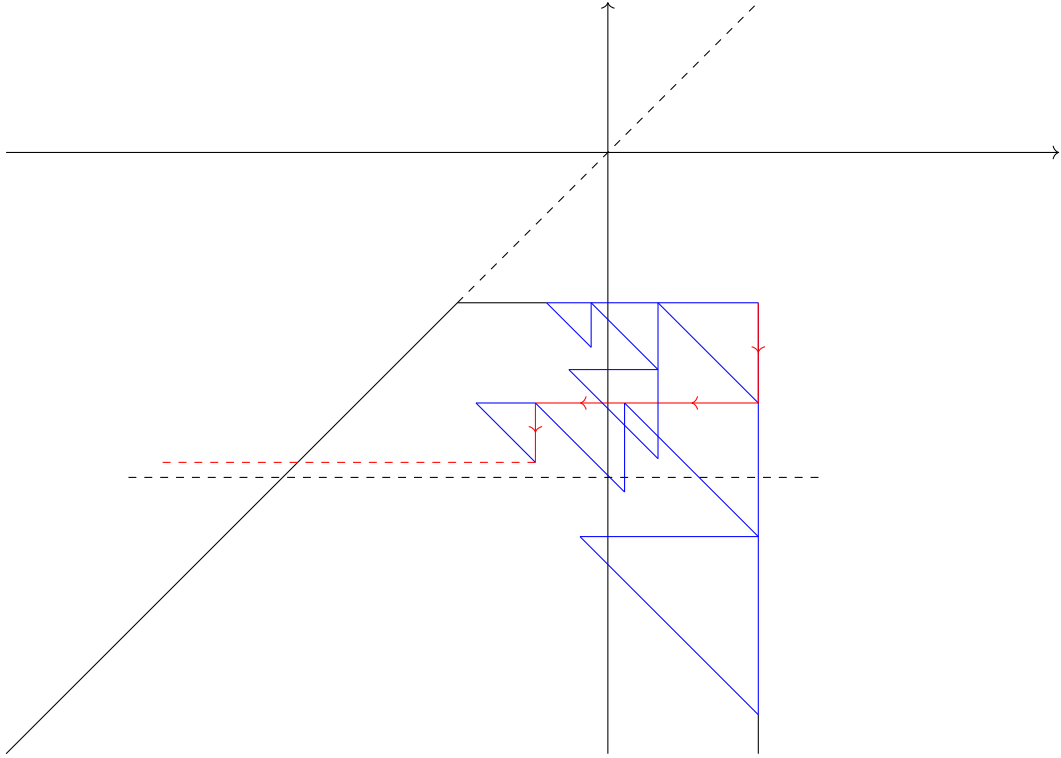
\begin{figure}[H]
    \centering
\begin{tikzpicture}[scale=2]
    % Draw axes

    \draw[->] (-4,0) -- (3,0) node[right] {}; % x-axis
    \draw[->] (0,-4) -- (0,1) node[above] {}; % y-axis

    % Draw the line y=x
    \draw[dashed] (-1,-1) -- (1,1); % Dashed section between -1 and 1
    \draw (-4,-4) -- (-1,-1);       % Solid section to the left of -1       % Solid section to the right of 1
\draw[blue] (1,-1) -- (-1/9-8/27,-1);
    \draw (-1/9-8/27,-1) -- (-1,-1);

    \draw [blue](1,-1) -- (1,-101/27);
    \draw (1,-4) -- (1,-101/27);
\draw [blue](1/3,-1) -- (1,-5/3);
\draw[red, dashed] (1/9-16/27,-1-2/3-32/81) -- (-3,-1-2/3-32/81);
\draw[dashed] (1.4,-1-2/3-32/81-1/10) -- (-3.2,-1-2/3-32/81-1/10);
\draw[red, postaction={decorate}, decoration={markings, mark=at position 0.5 with {\arrow{>}}}] 
        (1,-1) -- (1,-5/3) ;
        \draw[red, postaction={decorate}, decoration={markings, mark=at position 0.5 with {\arrow{>}}}] 
         (1,-5/3) --(1/9,-5/3) ;
\draw[blue] (1/9,-5/3) -- (1,-23/9);
\draw[red, postaction={decorate}, decoration={markings, mark=at position 0.5 with {\arrow{>}}}] 
        (1/9,-5/3) -- (1/9-16/27,-5/3) ;
\draw [blue](1-8/9-16/27-32/81,-5/3) -- (1/9-16/27,-5/3);
\draw [blue](1-8/9-16/27-32/81,-5/3) -- (1/9-16/27,-5/3-32/81);
    \draw[red, postaction={decorate}, decoration={markings, mark=at position 0.5 with {\arrow{>}}}] 
        (1-8/9-16/27,-5/3) -- (1/9-16/27,-1-2/3-32/81) ;
\draw[blue] (1/9,-5/3) -- (1/9,-5/3-16/27);
\draw[blue] (1/9-16/27,-5/3) -- (1/9,-5/3-16/27);
\draw[blue] (-5/27,-23/9) -- (1,-23/9);
\draw[blue] (-5/27,-23/9) -- (1,-101/27);

\draw[blue] (-1/9,-1) -- (1/3,-13/9);
\draw[blue] (1/3,-1) -- (1/3,-13/9);

\draw[blue] (-1/9,-1) -- (-1/9,-35/27);
\draw [blue](-1/9-8/27,-1) -- (-1/9,-35/27);

\draw[blue] (1/3-16/27,-13/9) -- (1/3,-13/9);
\draw [blue](1/3,-13/9-16/27) -- (1/3,-13/9);
\draw[blue] (1/3,-13/9-16/27) -- (1/3-16/27,-13/9);

\end{tikzpicture}\\

    \caption{Finding approximation, starting with $A=\emptyset$.}
    \label{fig:theorem}
\end{figure}

\begin{proof}
    Constructive. Let \(A\) be any finite sequence of \(2\)'s and \(1\)'s such that \(C^{-1}[A,(1)]>z\) (such \(A\) always exists, as \(z<-1\), and if we let \(A\) be the empty sequence, then \(C^{-1}[A,(1)]=-1\)). Look at the sequence of numbers:
    \begin{align}
        S_A= \Big(C^{-1}[A,(1)]\text{, }C^{-1}[A,2,(1)]\text{, }C^{-1}[A,2,2,(1)]\text{,...}\Big).
    \end{align}
    It's a strictly decreasing sequence diverging to \(-\infty\). Take the last element of it which is still bigger or equal to \(z\). Say it is \(C^{-1}[B,(1)]\), where \(B\) is the sequence \(A\) concatenated \(n\) times with the number \(2\). By construction \(C^{-1}[B,(1)]>z\), and \(C^{-1}[B,2,(1)]<z\). Let \(w\coloneqq C^{-1}[B,(1)]-z\). Since  \(C^{-1}[B,2,(1)]<z\), we get that \(w < C^{-1}[B,(1)]-C^{-1}[B,2,(1)]\).\\
    Consider the sequence of numbers:
    \begin{align}
        S_B= \Big(C^{-1}[B,(2)]\text{, }C^{-1}[B,1,(2)]\text{, }C^{-1}[B,1,1,(2)]\text{,...}\Big),
    \end{align}
    and the sequence of the differences of the consecutive terms of it:
    \begin{align}\label{diffs}
        S_B'=\Big(C^{-1}[B,(2)]-C^{-1}[B,1,(2)]\text{, }C^{-1}[B,1,(2)]-C^{-1}[B,1,1,(2)]\text{,...}\Big),
    \end{align}
    Let \(x \coloneqq C^{-1}[B,(2)]-C^{-1}[B,1,(2)]\). Apply statement \(3\) from Lemma \ref{statements}:
    \begin{align}
        C^{-1}[B,1,(2)]-C^{-1}[B,1,1,(2)]= 2/3 \cdot (C^{-1}[B,(2)]-C^{-1}[B,1,(2)]),
    \end{align}
to get that the sequence \(S_B'\) is of the form:
\begin{align}
   S_B'= \big(x\text{, }x\cdot(2/3)\text{, }x\cdot(2/3)^2\text{, }x\cdot(2/3)^3,...\big),
\end{align}
strictly decreasing and converging to \(0\). Take the first element of it which is smaller than \(w\). Let's say it's \(y = C^{-1}[B_1,(2)]-C^{-1}[B_1,1,(2)]\), where \(B_1\) is \(B\) concatenated with some number of \(1\)'s.\\

We claim that:
\begin{align}
    0<C^{-1}[B_1,2,(1)]-z<(1/3) \cdot (C^{-1}[B,(1)]-z) = (1/3)\cdot w
\end{align}
If this is true, we are done. We let \([B_1,2]\) be our sequence \(B\), we repeat what we've done and we find next sequence \([C_1,2]\) such that:
\begin{align}
    0<C^{-1}[C_1,2,(1)]-z&<(1/3) \cdot (C^{-1}[B_1,2,(1)]-z)\\
    &<(1/3^2) \cdot (C^{-1}[B,(1)]-z)\\
    &=(1/9)\cdot w,
\end{align}
and so on, until we get to a sequence \([Y]\) which is as close to \(z\) as we want.\\

So, let us prove our claim:\\

Since \(y\) was the first number in the sequence \(S_B'\) which was smaller than \(w\), we know that \((3/2)\cdot y >w\) (the previous term of this sequence is bigger that \(w\)).
Apply statement \(1\) from Lemma \ref{statements}:
    \begin{align}
        y= C^{-1}[B_1,(2)]-C^{-1}[B_1,1,(2)]=C^{-1}[B_1,(1)]-C^{-1}[B_1,2,(1)].
    \end{align}
But since \(B_1\) is  \(B\) concatenated with some number of \(1\)'s, we have that: \(C^{-1}[B_1,(1)]=C^{-1}[B,(1)]\). Hence:
\begin{align}
        C^{-1}[B,(1)]-C^{-1}[B_1,2,(1)] = y > (2/3)\cdot w.
    \end{align}
Therefore:
\begin{align}
        C^{-1}[B_1,2,(1)] -z <  C^{-1}[B,(1)]- (2/3)\cdot w-z
    \end{align}
But, by definition \(w = C^{-1}[B,(1)]-z\), so we finally get:
\begin{align}
        C^{-1}[B_1,2,(1)] -z <  (1/3)\cdot w.
    \end{align}
On the other hand, \(y\) is taken to be smaller than \(w\) and so:
\begin{align}
    C^{-1}[B,(1)]-C^{-1}[B_1,2,(1)] = y <w = C^{-1}[B,(1)]-z,
\end{align}
which means:
\begin{align}
    0< C^{-1}[B_1,2,(1)]-z,
\end{align}
which is what we wanted to prove.
    
\end{proof}

\section{Further Results}
\theorem We can meaningfully (but not uniquely) write each number \(z \leq -1\) as \(C[A]\) for \(A\) an infinite sequence consisting only of \(1\)'s and \(2\)'s.

\begin{proof}
    Follows from the constructive proof.
\end{proof}

\theorem Let \(A\) be any finite sequence of \(1\)'s and \(2\)'s. Then points \((1,-1)\), \(P(A)\), \(P(A,A)\), \(P(A,A,A)\), ..., all lie on the same line. Moreover, the point at which this line intersects the line \(y=x\) is precisely the point \(P((A))\).

\begin{proof}
    Calculations, but: also self similarity of the fractal structure of \(P_n\), for any \(n\).
\end{proof}

\remark This works even if the coordinates of \(P((A))\) are positive. I.e. even if the sequence of \(x\)- or \(y\)-coordinates of the points \(P(A)\), \(P(A,A)\), \(P(A,A,A)\), ... tends to \(-\infty\), we can still find the value of \(P((A))\). We do it by looking at the line through the points \((1,-1)\) and \(P(A)\) and calculating the point of intersection of it with the line \(y=x\).

\remark We proved that any eventually periodic sequence \(A\) defines a unique rational number \(C^{-1}(A)\) whose Collatz representation is precisely \(A\). This does not hold for all sequences which are not eventually periodic though. One can construct infinite sequences of positive integers (or even infinite sequences of only \(1\)'s and \(2\)'s) which are not equal to the Collatz representation of any number in \(\Q^{odd}\). It is similar to the case of continued fraction expansion of Laurent series with finitely many terms of positive degree, as in \cite{Malagoli}.

\remark Another nice result: given two finite sequences \(A\) and \(B\), we have that the numbers \(C^{-1}[(A,B)]\), \(C^{-1}[(B,A)]\) lie strictly between numbers \(C^{-1}[(A)]\) and \(C^{-1}[(B)]\) (unless the sequences \(A\) and \(B\) are equal, of course).

\remark In our notation, the Weak Collatz Conjecture (see \cite{tao}) can be reformulated as follows: the only positive integer in \(\Q^{odd}\) with absolutely periodic Collatz representation is \(1\). In general, it seems (see \cite{citation-key} for details) that the only integers in \(\Q^{odd}\) with absolutely periodic Collatz representations are \([1]\) (loop \(1\)), \([-1]\) (loop \(2\)), \([-5,-7]\) (loop \(3\)), \([-17,-25,-37,-55,-41,-61,-91]\) (loop \(4\)).
\end{document}